\documentclass[12pt]{article}
\usepackage{latexsym}
\usepackage{amssymb}
\usepackage{euscript}
\let\script\EuScript

\newtheorem{theorem}{Theorem}[section]
\newtheorem{definition}[theorem]{Definition}
\newtheorem{lemma}[theorem]{Lemma}

\newtheorem{notation}[theorem]{Notation}
\def\restrict{\upharpoonright}
\def\comp{\circ}
\def\N{\mathbb N}
\def\R{\mathbb R}

\def\Q{\mathbb Q}
\def\forces{\Vdash}
\def\nin{\notin}
\def\all{\forall}
\def\romu{{\mathbb R}^{\omega}/U}
\def\ro{\R^{\omega}}
\def\halmos{{}{\hfill $\Box$} \vskip 5pt \par}
\def\gaptwo{\langle \overrightarrow{\varphi}, \overrightarrow{\script{G}}, \overrightarrow{\theta}, \overrightarrow{\script{F}} \rangle}
\def\gapone{\langle \overrightarrow{\varphi}, \overrightarrow{\script{G}} \rangle}
\def\union{\cup}
\begin{document}
\setlength{\baselineskip}{21pt}
\title{Gap-2 Morass-definable $\eta_1$-orderings}
\author{Bob A. Dumas\\
University of Washington\\
Seattle, Washington 98195}
\date{December 9, 2016}

\bibliographystyle{plain}

\maketitle
\begin{abstract}
We prove that in the Cohen extension adding $\aleph_3$ generic
reals to a model of $ZFC+CH$ containing a simplified
$(\omega_1,2)$-morass, gap-2 morass-definable $\eta_1$-orderings
with cardinality $\aleph_3$ are order-isomorphic. Hence it is
consistent that $2^{\aleph_0}=\aleph_3$ and that morass-definable
$\eta_1$-orderings with cardinality of the continuum are
order-isomorphic. We prove that there are ultrapowers of $\R$ over
$\omega$ that are gap-2 morass-definable.  The constructions use a
simplified gap-2 morass, and commutativity with morass-maps and
morass-embeddings, to extend a transfinite back-and-forth
construction of order type $\omega_1$, to a function between
objects of cardinality $\aleph_3$.
\end{abstract}
\baselineskip = 18pt \setcounter{section}{0}
\section{Introduction}
An $\eta_1$-ordering without
endpoints, $\langle X,< \rangle$, is a linear-ordering for which
countably many consistent order constraints are
necessarily witnessed by an object in $X$.  That is, if $L\subset
X$ and $U\subset X$ are both countable, and
for every $l\in L$ and $u\in U$, $l<u$, then there is $x\in X$
such that for every $l\in L$ and $u\in U$, $l<x<u$.  We consider
only $\eta_1$-orderings without endpoints, that is $U$ or
$L$ above may be empty.

By the
Compactness Theorem, $\eta_1$-orderings exist in proliferation at
all infinite cardinalities, and there are many well-studied
examples of $\eta_1$-orderings that play a significant role in
logic, topology and analysis. It is an early result of Model
Theory that $\eta_1$-orderings having cardinality $\aleph_1$ are
order-isomorphic. This is proved with a classic back-and-forth
construction of the isomorphism.  The argument is by transfinite
construction of length $\omega_1$, requiring that only countably
many order commitments need be satisfied at any step of the
construction. Hence it is a consequence of the Continuum
Hypothesis ($CH$) that $\eta_1$-orderings having cardinality of
the continuum, $2^{\aleph_0}$, are order-isomorphic.

We are particularly interested in $\eta_1$-orderings without
endpoints bearing cardinality of the continuum. We seek to find
useful conditions on $\eta_1$-orderings in which $CH$ fails and
$\eta_1$-orderings bearing cardinality of the continuum are
order-isomorphic.  In [\ref{Dumas}] we showed that in a Cohen
extension adding $\aleph_2$-generic reals to a model of $ZFC +
CH$ containing a simplified $(\omega_1,1)$-morass, there is a
level order-isomorphism between morass-definable
$\eta_1$-orderings with cardinality of the continuum (Theorem 5.4,
[\ref{Dumas}]).  The simplified-morass plays a critical role in
that we reduce the element-wise construction of a function between
orderings having cardinality $\aleph_2$ to a construction of
length $\omega_1$.

In this paper, working in the Cohen extension adding generic reals
indexed by $\omega_3$ to a model of $ZFC + CH$ containing a
simplified $(\omega_1,2)$-morass, we define a level
order-isomorphism between gap-2 morass-definable
$\eta_1$-orderings having cardinality $2^{\aleph_0}=\aleph_3$.
Our strategy is as follows.  Given a simplified gap-2 morass,
$\gaptwo$, we construct a sequence of level order-preserving bijections,
$\langle F_{\beta}\mid \beta<\omega_1 \rangle$, on
the ``fake" morasses below $\omega_1$ in a manner
similar to the construction of [\ref{Dumas}], but with the
additional provision that the construction is closed under the embeddings of $\script{F}_{\beta \gamma}$,
for all $\beta<\gamma<\omega_1$.
Then $F=\bigcup_{\beta<\omega_1} F_{\beta}$ is a level
term function that is forced to be an order-isomorphism
between those elements of the $\eta_1$-orderings
that are also in the generic extension adding reals indexed by
$\omega_1$. Furthermore, we are able to use the embeddings of $\script{F}_{\beta \omega_1}$,
$\beta<\omega_1$,
to extend $F$ to a level term function that is forced to
be an order-isomorphism between
$\eta_1$-orderings having cardinality of the continuum in a model
of set theory with $\aleph_3$-generic reals.

In the next paper we show that morass-definable and gap-2
morass-definable $\eta_1$-ordered real-closed fields bearing
cardinality of the continuum are isomorphic in the Cohen
extensions adding $\aleph_2$ or $\aleph_3$ generic reals,
by an $\R$-linear order-preserving isomorphism. The
role of $\eta_1$-ordered real-closed fields in the subject of
automatic continuity (the existence of discontinuous homomorphisms
of $C(X)$, the algebra of continuous real-valued functions on an
infinite compact Hausdorff space, $X$) has been well-explored in
$[\ref{Dales1}]$, $[\ref{Dales2}]$, $[\ref{Esterle1}]$,
$[\ref{Esterle2}]$, $[\ref{Esterle3}]$ and $[\ref{Woodin}]$.  We
use the techniques of this paper to show that it is consistent
that the continuum has cardinality $\aleph_3$, and that there
exists a discontinuous homomorphism of $C(X)$, for any infinite
compact Hausdorff space $X$.
\section{The Simplified Gap-2 Morass}\label{sec3}
For $\kappa$ a regular cardinal, we define a simplified
$(\kappa,2)$-morass as in Definition 1.3 [\ref{Velleman2}].
\begin{definition} $(Simplified \: \: (\kappa,2)-morass)$ \label{def3.1}
The structure $\gaptwo$ is a simplified $(\kappa,2)$-morass
provided it has the following properties:
\begin{enumerate}
\item $\gapone$ is a neat simplified $(\kappa^+,1)$-morass
$[\ref{Velleman1}]$. \item $\all \alpha<\beta\leq \kappa$,
$\script{F}_{\alpha \beta}$ is a family of embeddings (see page
172, $[\ref{Velleman2}]$) from $\langle \langle
\varphi_{\zeta}\mid \zeta<\theta_{\alpha}\rangle, \langle
\script{G}_{\zeta \xi} \mid \zeta<\xi\leq \theta_{\alpha}\rangle
\rangle$ to $\langle \langle \varphi_{\zeta}\mid
\zeta<\theta_{\beta}\rangle, \langle \script{G}_{\zeta \xi} \mid
\zeta<\xi\leq \theta_{\beta}\rangle \rangle$. \item $\all
\alpha<\beta<\kappa \: \: (\mid \script{F}_{\alpha \beta}
\mid<\kappa)$.
\item $\all \alpha<\beta<\gamma\leq \kappa \: \: (\script{F}_{\alpha \gamma}=\{ f \comp g\mid f\in \script{F}_{\beta \gamma}, g\in \script{F}_{\alpha \beta} \} )$. Here $f\comp g$ is defined by:\\
    \[ (f\comp g)_{\zeta}=f_{g(\zeta)} \comp g_{\zeta} \; \; \; \; for \; \zeta\leq \theta_{\alpha}, \]
    \[ (f\comp g)_{\zeta \xi}=f_{g(\zeta) g(\xi)} \comp g_{\zeta \xi} \; \; \; \; for \; \zeta<\xi\leq \theta_{\alpha}. \]
\item $\all \alpha<\kappa$, $\script{F}_{\alpha \alpha+1}$ is an
amalgamation (see page 173 $[\ref{Velleman2}]$). \item If
$\beta_1, \beta_2 <\alpha \leq \kappa$, $\alpha$ a limit ordinal,
$f_1\in \script{F}_{\beta_1 \alpha}$ and $f_2\in
\script{F}_{\beta_2, \alpha}$, then $\exists \beta (\beta_1,
\beta_2 < \beta <\alpha$ and $\exists f_1'\in \script{F}_{\beta_1
\beta} \: \exists f_2'\in \script{F}_{\beta_2 \beta} \: \exists
g\in \script{F}_{\gamma \alpha} (f_1=g\comp f_1'$ and $f_2=g\comp
f_2'))$. \item If $\alpha\leq \kappa$ and $\alpha$ is a limit
ordinal, then:
\begin{enumerate}
\item $\theta_{\alpha}=\bigcup \{ f[\theta_{\beta}] \mid
\beta<\alpha$, $f\in \script{F}_{\beta \alpha} \}$. \item $ \all
\zeta \leq \theta_{\alpha}$, $\varphi_{\zeta}=\bigcup \{
f_{\bar{\zeta }}[\varphi_{\bar{\zeta}}] \mid \exists \beta<\alpha
(f\in \script{F}_{\beta \alpha}$, $f(\bar{\zeta})=\zeta) \}$.
\item $\all \zeta<\xi \leq \theta_{\alpha}$, $\script{G}_{\zeta
\xi}=\bigcup \{ f_{\bar{\zeta} \bar{\xi}}[\script{G}_{\bar{\zeta}
\bar{\xi}}] \mid \exists \beta<\alpha \: (f\in \script{F}_{\beta
\alpha}$, $f(\bar{\zeta})=\zeta$, $f(\bar{\xi})=\xi ) \}$.
\end{enumerate}
\end{enumerate}
\end{definition}
We consider the case in which $\kappa=\omega_1$.
If $\alpha<\beta \leq \omega_1$, and $f\in \script{F}_{\alpha \beta}$,
then following the notational simplification of Velleman, we consider
$f$ to be a triple:
\[ \langle f, \{ f_{\zeta} \mid \zeta<\theta_{\alpha} \},
\{ f_{\zeta \xi} \mid \zeta<\xi\leq \theta_{\alpha} \} \rangle. \]
In the triple above, $f:\theta_{\alpha}+1 \to \theta_{\beta}+1$
is an order-preserving injection with $f(\theta_{\alpha})=
\theta_{\beta}$.  We will refer to this as the first component of the embedding.
For $\zeta \leq \theta_{\alpha}$, $f_{\zeta}:\varphi_{\zeta} \to \varphi_{f(\zeta)}$ is
an order-preserving injection.  We refer to this as the second component
of the embedding (corresponding to $\zeta$).  Finally for $\zeta<\xi\leq
\theta_{\alpha}$, $f_{\zeta \xi}:\script{G}_{\zeta \xi} \to
\script{G}_{f(\zeta) f(\xi)}$ is a function between morass maps
of $\gapone$.  We refer to this as the third component of the embedding
(corresponding to $\zeta$ and $\xi$).  Embeddings satisfy a number of
regularity and commutativity conditions that make them a
practical tool for extending our
$\omega_1$-inductive construction to a construction
of cardinality $\aleph_3$.

\section{Morass Maps and Morass-Embeddings}
Let $M$ be a model of $ZFC + CH$ containing a simplified
$(\omega_1,2)$-morass, $\gaptwo$. Let $P$ be the poset
$Fn(\omega_3 \times \omega, 2)$, $G$ be $P$-generic over $M$ and
$M^P$ be the forcing language of the poset $P$ in $M$.
The construction of this paper is dependent on details of
the indexing set of $P$.  We are required to construct
a function between terms in the forcing language that
is sensitive to the precise subset of $\omega_3$ required for
the construction of carefully selected terms of $M^P$.
In general, if $S\subset \omega_3$, we let $P(S)$
be the poset $Fn(S\times \omega, 2)$.  In order to easily
associate these partial languages with constructions along the
morass, for $\beta\leq \omega_3$, let
$P_{\beta}=P(\varphi_{\theta_{\beta}}\times \omega,2)$.
If $G$ is $P$-generic over $M$, let $G(S)$ be the factor of
$G$ that is $P(S)$-generic over $M$ (similarly for $G_{\beta}$).

We will construct a function between sets of terms with
cardinality $\aleph_3$ in the forcing language, $M^P$, by applying
the second components of embeddings of a simplified gap-2 morass to a function on a set
of terms in $M^{P_{\omega_1}}$ with cardinality $\aleph_1$.  If
$\gaptwo$ is a simplified $(\omega_1, 2)$-morass, then $\gapone$
is a simplified $(\omega_2,1)$-morass. It is easy to verify that
every countable subset of $\omega_3$ is in the image of a single morass
function of $\gapone$.  If every
countable subset of $\omega_3$ were in the range of a
morass map from a countable vertex, then the construction of
Theorem 5.4 [\ref{Dumas}] would suffice to prove that
morass-definable $\eta_1$-orderings are order-isomorphic in $M[G]$.
However, not every countable subset of $\omega_3$
is anticipated by a morass map from
a countable vertex of $\gapone$. Instead we employ the embeddings of
a simplified $(\omega_1,2)$-morass to anticipate countable subsets of
$\omega_3$ by countable subsets of $\omega_1$, and thereby
construct a function between sets of terms in the forcing language
having cardinality $\aleph_3$ by an inductive construction of
length $\omega_1$.

The basic strategy of the central construction of this paper is to define
a sequence of functions on terms in the forcing language adding generic reals
indexed by $\omega_1$ with the enhanced back-and-forth argument
similar to the argument in [\ref{Dumas}]. The term functions so constructed will be
level and satisfy certain closure conditions with the second components
of morass-embeddings.
We use the embeddings of $\script{F}_{\beta \omega_1}$ to ``lift" the
construction to a function on terms in the forcing language adding
generic reals indexed by $\omega_3$.  We use the following pair of Lemmas
due to Velleman.
\begin{lemma}\label{lem3.2} (Velleman). Let $\alpha<\beta\leq \omega_2$, $g_1, g_2\in
\script{G}_{\alpha \beta}$, $\tau_1, \tau_2\in \varphi_{\alpha}$ and
$g_1(\tau_1)=g_2(\tau_2)$.  Then $\tau_1=\tau_2$ and $g_1\restrict_{\tau_1+1}
=g_2\restrict_{\tau_2+1}$.
\end{lemma}
The gap-2 version of this is Lemma 2.2 [\ref{Velleman2}]:
\begin{lemma}(Velleman)\label{lem3.3}
Suppose $\alpha<\beta\leq \kappa$, $f_1,f_2\in \script{F}_{\alpha \beta}$,
$\zeta_1, \zeta_2\leq \theta_{\alpha}$, and $ssup(f[\zeta_1])=ssup(f[\zeta_2])$.
Then $\zeta_1=\zeta_2$, $f_1\restrict_{\zeta_1}=f_2\restrict_{\zeta_2}$ and
$(f_1)_{\xi}=(f_2)_{\xi}$ for all $\xi<\zeta_1$.
\end{lemma}
We show that every countable subset of $\omega_3$ is in the
range of a morass function, albeit not necessarily from
a countable vertex.
\begin{lemma} Let $T$ be a countable subset of $\omega_3$.  Then there is $\zeta<\omega_2$
and $g\in \script{G}_{\zeta \omega_2}$ such that $T\subseteq g[\varphi_{\zeta}]$.
\end{lemma}
Proof.
Let $\tau=sup(T)$.  We assume that $\tau\in T$.  Enumerate the elements of
$T$, $\langle \tau_n \mid n\in \omega\rangle$, with $\tau=\tau_0$.  By condition P5 of the definition
of a simplified gap-1 morass
there is $\zeta_0<\omega_2$, $\sigma_0\in \varphi_{\zeta_0}$ and $g_0\in \script{G}_{\zeta_0 \omega_2}$
such that $g_0(\sigma_0)=\tau$.  For $n\in \omega$, there is $\zeta'<\omega_2$,
$\sigma'\in \varphi_{\zeta'}$ and $g'\in \script{G}_{\zeta' \omega_2}$ such that
$g'(\sigma')=\tau_n$.  By condition P4 of the definition of a simplified gap-1 morass
there is $\zeta_n \geq max(\zeta_0,\zeta')$, $h_0\in \script{G}_{\zeta_0 \zeta_n}$,
$h_n'\in \script{G}_{\zeta' \zeta_n}$ \ and $g_n\in \script{G}_{\zeta_n \omega_2}$
such that $g_0=g_n \comp h_0$ and $g'=g_n \comp h_n'$.
We observe that
\[ g_n(h_0(\sigma_0))=\tau \]
and
\[ g_n(h_n'(\sigma_n))=\tau_n. \]
Hence if $n\in \omega$, $\tau$ and $\tau_n$ are in the image
of a single morass map
$g_n\in \script{G}_{\zeta_n \omega_2}$.

Let $\omega_2>\zeta>\zeta_n$ for all $n\in \omega$.
By condition P2 of the simplified gap-1 morass, there is $g\in \script{G}_{\zeta \omega_2}$
and $g_0^*\in \script{G}_{\zeta_0 \zeta}$ such that
\[ g_0=g\comp g_0^*. \]
We observe that
\[ g_0(h_0(\sigma_0))=g(g_0^*(h_0(\sigma_0)))=\tau. \]
Let $n\in \N$, $n>0$, $g_n^*\in \script{G}_{\zeta_n \zeta}$
and $g^*\in \script{G}_{\zeta \omega_2}$ be such that
\[ g_n=g^*\comp g_n^*. \]
Then
\[ g^*(g_n^*(h_0(\sigma_0)))=\tau. \]
So by the previous lemma,
\[ g_n^*(h_0(\sigma_0))=g_0^*(h_0(\sigma_0)) \]
and
\[ g\restrict_{g_0^*(h_0(\sigma_0))+1}=g^*\restrict_{g_0^*(h_0(\sigma_0))+1}. \]
In particular, $\tau_n<\tau$ and $g\comp g^*_n(\sigma_n)=\tau_n$.
Hence, for all $n\in \omega$,
\[ \tau_n\in g[\varphi_{\zeta}]. \]
\halmos
\begin{lemma}
Assume
\begin{enumerate}
\item $\zeta<\omega_2$ and $g\in \script{G}_{\zeta \omega_2}$
\item $\beta_1, \beta_2<\omega_1$
\item For $i=1,2$, $f_i\in \script{F}_{\beta_i \omega_1}$
\item For $i=1,2$, $\bar{\zeta_i}<\theta_{\beta_i}$ and $f_{i}(\bar{\zeta_i})=\zeta$.
\item For $i=1,2$, $\bar{g}_i\in \script{G}_{\bar{\zeta_i} \theta_{\beta_i}}$ and
$\script{G}_{\bar{\zeta_i} \theta_{\beta_i}}(\bar{g}_i)=g$.
\end{enumerate}
Then there is countable $\beta\geq \beta_i$, $i=1,2$,
$f\in \script{F}_{\beta \omega_1}$, $\bar{\zeta}\in \theta_{\beta}$,
and $\bar{g}\in \script{G}_{\bar{\zeta} \theta_{\beta}}$ such that
\begin{enumerate}
\item $f(\bar{\zeta})=\zeta$
\item $f_{\bar{\zeta} \theta_{\beta}}(\bar{g})=g$
\item For $i=1,2$, there are $f_i'\in \script{F}_{\beta_i \beta}$
such that $f_i=f\comp f_i'$.
\end{enumerate}
\end{lemma}
Proof:
Let $\zeta$, $g$, $\beta_i$, $f_i$, $\bar{\zeta}_i$ and $\bar{g}_i$
($i=1,2$), satisfy the hypotheses.  By Condition 6
of Definition \ref{def3.1}, there are $\beta$, $f$ and $f_i'$ ($i=1,2$)
such that $f_i=f\comp f_i'$.  Let $\bar{\zeta}=f_1'(\bar{\zeta}_1)$.
Then $\bar{\zeta}<\theta_{\beta}$ and
\[ f(\bar{\zeta})=f\comp f_1'(\bar{\zeta}_1)=f_1(\bar{\zeta})=\zeta. \]
Since $f$ is an injection and $f_2(\bar{\zeta}_2)=f\comp f_2'(\bar{\zeta}_2)=\zeta$,
\[ f_2'(\bar{\zeta}_2)=\bar{\zeta}. \]
Let $\bar{g}=(f_1')_{\bar{\zeta}_1 \bar{\zeta}}(\bar{g}_1)$.
Then
\[ f_{\bar{\zeta} \theta_{\beta}}(\bar{g})=f_{\bar{\zeta} \theta_{\beta}}((f_1')_{\bar{\zeta}_1 \bar{\zeta}}(\bar{g}_1))
=f_1(\bar{g}_1)=g. \]
Since $f$ is an embedding,
\[ \bar{g}=(f_2')_{\bar{\zeta}_2 \bar{\zeta}}(\bar{g}_2). \]
\halmos
\begin{definition}(Continuation)
Let $\alpha<\beta<\gamma\leq\omega_1$ and $f\in \script{F}_{\alpha \gamma}$.
If $f'\in \script{F}_{\alpha \beta}$ and $f^*\in \script{F}_{\beta \gamma}$
are such that $f=f^*\comp f'$, then we say that $f^*$ is the
$\beta$-continuation of $f$.
\end{definition}
Condition 4 of Definition \ref{def3.1} states that,
for any $\alpha<\beta<\gamma\leq\omega_1$, if
$f\in \script{F}_{\alpha \gamma}$, then $f$ has
a $\beta$-continuation.  Condition 6 of Definition \ref{def3.1}
states that if $\beta_1,\beta_2<\gamma\leq \omega_1$,
$\gamma$ a limit ordinal,
$f_1\in \script{F}_{\beta_1 \gamma}$ and $f_2\in \script{F}_{\beta_2 \gamma}$
then there is $\beta \geq \beta_1, \beta_2$
such that $f_1$ and $f_2$ have a common $\beta$-continuation.
\begin{lemma}\label{lem3.6}
Let $\zeta<\omega_2$, $g\in \script{G}_{\zeta \omega_2}$ and
$\langle (\beta_n', f'_n, \zeta_n', g_n')\mid n\in \N \rangle$ be a sequence
satisfying:
\begin{enumerate}
\item $f'_n\in \script{F}_{\beta_n' \omega_1}$
\item For all $n$, $\zeta_n'<\theta_{\beta_n'}$
and $f'_n(\zeta_n')=\zeta$
\item For all $n$, $g_n'\in \script{G}_{\zeta_n' \theta_{\beta_n'}}$
and $(f'_n)_{\zeta_n' \theta_{\beta_n'}}(g_n')=g$.
\end{enumerate}
Then there is a sequence $\langle (\beta_n, f_n, \zeta_n, g_n) \mid n\in \N \rangle$
satisfying:
\begin{enumerate}
\item For all $m<n$, $\beta_n\geq \beta_n'$ and $\beta_n\geq \beta_m$
\item For all $n$, $f_n\in \script{F}_{\beta_n \omega_1}$
\item For all $n$, $f_n(\zeta_n)=\zeta$
\item For all $n$, $(f_n)_{\zeta_n \theta_{\beta_n}}(g_n)=g$
\item For all $n$, $f_n$ is the $\beta_n$-continuation of $f'_n$
\item For all $m<n$, $f_n$ is the $\beta_n$-continuation of $f_m$.
\end{enumerate}
The sequence $\langle (\beta_n, f_n, \zeta_n, g_n) \mid n\in \N \rangle$ is
called an embedding sequence of $g$.
\end{lemma}
Proof:
Let $\langle (\beta_n', f'_n, \zeta_n', g_n')\mid n\in \N \rangle$ satisfy
the hypothesis of the Lemma.
We construct $\langle (\beta_n, f_n, \zeta_n, g_n) \rangle$ by induction.
Let $(\beta_0', f'_0, \zeta_0', g_0')=(\beta_0, f_0, \zeta_0, g_0)$.
Assume $N\in \N$ and that $\langle (\beta_n, f_n, \zeta_n, g_n)\mid n\leq N \rangle$
satisfies the conclusion of the Lemma beneath $N$.  By Condition 4
and Condition 6 of Definition \ref{def3.1}, there is countable $\beta_{N+1}>\beta_N, \beta'_{N+1}$
and $f_{N+1}\in \script{F}_{\beta_{N+1} \omega_1}$ that is the $\beta_{N+1}$-continuation
of $f_N$ and $f'_{N+1}$.  Let $f_N=f_{N+1}\comp f'$ where $f'\in \script{F}_{\beta_N \beta_{N+1}}$
and $f^*\in \script{F}_{\beta_{N+1} \omega_1}$.  Let $\zeta_{N+1}=f'(\zeta_N)$.  Then
\[ f_{N+1}(\zeta_{N+1})=f_{N+1}\comp f'(\zeta_N)=f_N(\zeta_N)=\zeta. \]
Let $g_{N+1}=f'_{\zeta_N \zeta_{N+1}}(g_N)$.  Then
\[ (f_{N+1})_{\zeta_{N+1} \theta_{\beta_N}}(g_{N+1})=
(f_{N+1})_{\zeta_{N+1} \theta_{\beta_{N+1}}}(f'_{\zeta_N \zeta_{N+1}}(g_N))=
(f_N)_{\zeta_N \theta_{\beta_N}}(g_N)=g. \]
For any $n<N$, $F_{N+1}$ is the $\beta_{N+1}$-continuation
of $f_n$.
\halmos
\begin{definition}(Complete embedding sequence)
Let $\zeta<\omega_2$ and $g\in \script{G}_{\zeta \omega_2}$
and $S=\langle (\beta_n, f_n, \zeta_n, g_n) \mid n\in \N \rangle$
be an embedding sequence of $g$ and $D\subseteq \varphi_{\zeta}$.
Then $S$ is a complete embedding
sequence for $D$, with respect to $g$, provided
that $D\subseteq \bigcup_{n\in \N} (f_n)_{\zeta_n}[\varphi_{\zeta_n}]$.
\end{definition}
If $S$ is a complete embedding sequence of $g$, then each
embedding $f_n$ gives partial information about $g$.  That is
\[ g\comp (f_n)_{\zeta_n}=(f_n)_{\theta_{\beta_n}}\comp g_n. \]
Later embeddings in the sequence, $f_m$ ($m>n$) necessarily agree with earlier embeddings
about $g$, on $(f_n)_{\zeta_n}[\varphi_{\zeta_n}]$, but in general
will provide information about $g$ on a larger subset of $\varphi_{\zeta}$,
$(f_m)_{\zeta_m}[\varphi_{\zeta_m}]$.  A complete embedding sequence will
provide full information on $D$.

The second term of the embeddings
also provides information about $g[D]$.
For all $n \in \N$,
\[ (f_n)_{\theta_{\beta_n}}\comp g_n
=g\comp (f_n)_{\zeta_n}. \]
If $m>n$, then
\[ (f_n)_{\zeta_n}[\varphi_{\zeta_n}]\subseteq
(f_m)_{\zeta_m}[\varphi_{\zeta_m}] \]
and
\[ (f_n)_{\theta_{\beta_n}}\comp g_n[\varphi_{\zeta_n}]
\subseteq (f_m)_{\theta_{\beta_m}}\comp g_m[\varphi_{\zeta_m}]. \]
\begin{lemma}\label{lem3.9}
Let $T\subseteq \omega_3$ be countable,
$\zeta<\omega_2$,
$g\in \script{G}_{\zeta \omega_2}$
and $D\subseteq \varphi_{\zeta}$ such
that $g[D]=T$.
Then there is
$\gamma<\omega_1$, $\bar{\zeta}<\theta_{\gamma}$,
$h\in \script{F}_{\gamma \omega_1}$ and
$\bar{g}\in \script{G}_{\bar{\zeta} \theta_{\gamma}}$
such that
\[ D\subseteq h_{\bar{\zeta}}[\varphi_{\bar{\zeta}}] \]
\[ T\subseteq h_{\theta_{\gamma}}[\varphi_{\theta_{\gamma}}] \]
and
\[ h_{\bar{\zeta} \theta_{\gamma}}(\bar{g})=g. \]
\end{lemma}
Proof:
Let $T\subseteq \omega_3$, $\zeta<\omega_2$,
$D\subseteq \varphi_{\zeta}$ and $g\in \script{G}_{\zeta \omega_2}$
satisfy the hypothesis of the lemma.
By Lemma \ref{lem3.6}
there is a complete embedding sequence for
$D$ with respect to $g$, $S=\langle (\beta_n, f_n, \zeta_n, g_n) \mid n\in \N \rangle$.
For all $m>n\geq 0$, $f_m$ is
the $\beta_m$-continuation of $f_n$.
Therefore, for all $n\in \N$, there is $\rho_n\in \varphi_{\zeta_n}$ such that
\[ (f_0)_{\varphi_{\zeta_0}}(\rho_0)=
(f_n)_{\varphi_{\zeta_n}}(\rho_n)=\sigma. \]
Let $\gamma$ be a countable and, for all $n\in \N$, $\beta_n<\gamma$.
For $n\in \N$, let $h_n\in \script{F}_{\gamma \omega_1}$
be the $\gamma$-continuation of $f_n$.
Let $h=h_0$.
For each $n\in \N$,
there is $f_n^*\in \script{F}_{\beta_n \gamma}$
such that
\[ f_n=h_n\comp f_n^*. \]
Let $\bar{\zeta}=f_0^*(\zeta_0)$.
For all $n>0$, $f_n$ is the
$\beta_n$-continuation of $f_0$.  Hence
\[ f_n^*(\zeta_n)=\bar{\zeta}. \]
We observe that
\[ ssup(h[\bar{\zeta}+1])=\zeta=ssup(h_n[\bar{\zeta}+1]). \]
By a Lemma \ref{lem3.3},
\[ h\restrict_{\bar{\zeta}+1}=h_n\restrict _{\bar{\zeta}+1} \]
and
\[ h_{\bar{\zeta}}=(h_n)_{\bar{\zeta}}. \]
Since $S$ is a complete embedding sequence
with respect to $g$,
\[ D\subseteq \bigcup_{n\in \N} (h_n)_{\bar{\zeta}}[\varphi_{\bar{\zeta}}]. \]
Then
\[ D\subseteq h_{\bar{\zeta}}[\varphi_{\bar{\zeta}}]. \]
Let $\bar{g}=(f^*_0)_{\zeta_0 \theta_0}(g_0)\in \script{G}_{\bar{\zeta} \theta_{\gamma}}$.
Then
\[ h_{\bar{\zeta} \theta_{\gamma}}(\bar{g})=
(f_0)_{\zeta_0 \beta_0}(g_0)=g. \]
So
\[ h_{\theta_{\gamma}}\comp \bar{g}[\varphi_{\bar{\zeta}}]
=g\comp h_{\bar{\zeta}}[\varphi_{\bar{\zeta}}]. \]
However
\[ \bar{g}[\varphi_{\bar{\zeta}}]\subseteq \varphi_{\theta_{\gamma}}. \]
Therefore
\[ T\subseteq h_{\theta_{\gamma}}[\varphi_{\theta_{\gamma}}]. \]
\halmos
\begin{definition}(Compatible Maps)\label{def3.10}
Let $\theta \leq \kappa$ be ordinals, $f:\theta \to \kappa$
and $g:\theta \to \kappa$ be injections.  The ordinal maps $f$ and $g$ are
compatible provided that for any $\alpha, \beta \in \theta$,
$f(\alpha)=g(\beta)$ implies $\alpha=\beta$.
\end{definition}
If $f$ and $g$ are compatible, then $(f\union g)^{-1}$ is a
well-defined function.
Lemma \ref{lem3.2} implies that morass-maps with the same domain and codomain
are compatible.
\begin{lemma}\label{lem3.12}
Let $\beta<\omega_1$, $\eta$ be the splitting point of the right-branching embedding of
$\script{F}_{\beta \beta+1}$, $f,g\in \script{F}_{\beta \beta+1}$
and $\bar{f},\bar{g}\in \script{G}_{\eta \theta_{\beta}}$.
Then $f_{\theta_{\beta}}\comp \bar{f}$ and $g_{\theta_{\beta}} \comp \bar{g}$ are compatible.
If $p\in P_{\beta}$ then $f_{\theta_{\beta}}(p)$ and $g_{\theta_{\beta}}(p)$ are
compatible conditions in $P$.
\end{lemma}
Proof.
If $f$ and $g$ are left-branching embeddings then
$f_{\theta_{\beta}}\comp \bar{f}, g_{\theta_{\beta}}\comp \bar{g}\in \script{G}_{\eta \theta_{\beta+1}}$.
By Lemma \ref{lem3.2} they are compatible.  Assume that $f$ is the right-branching embedding
of $\script{F}_{\beta \beta+1}$.  Then
\[ f_{\theta_{\beta}}\comp \bar{f}=f_{\eta \theta_{\beta}}(\bar{f})\comp f_{\eta}\in \script{G}_{\eta \theta_{\beta+1}}. \]
Therefore, by Lemma \ref{lem3.2}, $f_{\theta_{\beta}}\comp \bar{f}$ and
$g_{\theta_{\beta}}\comp \bar{g}$ are compatible.

Let $p\in P_{\beta}$, $g$ and $h$ be second components of embeddings of $\script{F}_{\beta \beta+1}$.
So $g$ and $h$ are members of $\script{G}_{\theta_{\beta} \theta_{\beta+1}}$
or are equal to $f_{\theta_{\beta}}$.
Let $(\alpha,k,s)\in g(p)$, $(\alpha,k,t)\in h(p)$,
$g(\alpha')=\alpha$ and $h(\alpha^*)=\alpha$.
Then
\[ (\alpha',k,s)\in p \]
and
\[ (\alpha^*,k,t)\in p. \]
Morass maps are compatible, hence if $g,h\in \script{G}_{\theta_{\beta} \theta_{\beta+1}}$
then $\alpha'=\alpha^*$ and $s=t$.

Suppose that $h=f_{\theta_{\beta}}$.  Then there is
$\bar{\alpha}\in \varphi_{\eta}$ (where $\eta$ is the splitting point of $f$)
and $\bar{h}\in \script{G}_{\eta \theta_{\beta}}$
such that $\bar{h}(\bar{\alpha})=\alpha^*$.
So
\[ f_{\theta_{\beta}}\comp \bar{h}(\bar{\alpha})=
f_{\eta \theta_{\beta}}(\bar{h})\comp f_{\eta}(\bar{\alpha})=\alpha. \]
Since $f_{\eta \theta_{\beta}}(\bar{h}) \in \script{G}_{\theta_{\beta} \theta_{\beta+1}}$,
\[ f_{\eta}(\bar{\alpha})=\alpha'. \]
Therefore $(\alpha',k,t)\in p$ and $s=t$.
\halmos
\section{The Gap-1 Construction}
In [\ref{Dumas}], we constructed an
order-isomorphism between morass-definable
$\eta_1$-orderings in the generic extension
adding $\aleph_2$-generic reals
of a model of $ZFC+CH$ containing a
simplified $(\omega_1,1)$-morass.
This was accomplished by an
inductive construction of length
$\omega_1$ on terms in the forcing language
for adding $\aleph_1$ generic reals, subject
to numerous technical constraints, and using
the morass maps to ``lift" the construction
to the generic extension adding $\aleph_2$
generic reals.  The benefit of the gap-1 morass
is that the necessary back-and-forth construction
occurs strictly in the forcing language for adding only
$\aleph_1$ generic reals, relying on the morass-maps
for the completion of the overall construction.

The conditions required for the gap-1 construction were that the
$\eta_1$-orderings be morass-definable.
We present some adjustments of
definitions from the gap-1 construction
that will serve us in the gap-2 construction.

Let $\gaptwo$ be a simplified $(\kappa, 2)$-morass.
\begin{notation}($P_{\zeta}$, $P(A)$, $M^{P_{\zeta}}$, $M^{P(A)}$)
If $\zeta \leq \kappa$, $P_{\zeta}$ is the poset adding generic reals
indexed by $\varphi_{\theta_{\zeta}}$ and $M^{P_{\zeta}}$ is the
set of terms in the forcing language of $P_{\zeta}$.
For $A \subseteq \kappa^{++}$, $P(A)$ is the poset
for adding generic reals indexed by $A$, and
$M^{P(A)}$ is the set of terms in the forcing language
of $P(A)$.
\end{notation}
So $P=P_{\kappa}=P(\kappa^{++})$.
\begin{definition}(Support)
Let $\tau\in M^P$.  The support of $\tau$, $supp(\tau)$, is the minimal
subset of $\omega_3$, $A$, such that $\tau\in M^{P(A)}$.
We say that $\tau$ has countable support provided that $A$ is countable.
\end{definition}
Let $A\subseteq B\subseteq \omega_3$. We consider
$P(B)=P(A)\times P(B\setminus A)$.  If $G$ is
$P$-generic over $M$, then $M[G]=M[G(A)\times H]$
where $G(A)$ is $P(A)$-generic over $M$ and
$H$ is $P(B\setminus A)$-generic over $M[G(A)]$.
If $\tau$ is a term of the forcing language of $P(A)$,
the value of $\tau$ in the generic
extensions will be an element of $M[G(A)]$.
\begin{definition}(Strict Support)
Let $B\subseteq \kappa$ and $\tau\in M^{P(B)}$.  The term
$\tau$ has strict support $B$ if there is no proper
subset $A\subset B$, $p\in P$ and $\sigma \in M^{P(A)}$ such that
$p\forces \sigma=\tau$.
\end{definition}
If a term, $\tau$, has strict support $B$, then $val_{G(B)}(\tau)$
will not be an element of $M[G(A)]$, for any proper
subset $A$ of $B$.
\begin{definition}(Discerning Set of Terms)
A set of terms, $X\in M$, $X\subseteq M^P$, is discerning provided that every
term of $X$ has strict support.
\end{definition}
\begin{definition}(Level Function)
Let $X$ and $Y$ be discerning sets of terms and $\phi:X\to Y$.  The function
$\phi$ is level if for any $x\in X$, $x$ and $\phi(x)$ have identical strict support.
\end{definition}
\begin{definition}(Morass-Commutative)
Let $\kappa$ be regular, $\gapone$ be a simplified
$(\kappa,1)$-morass, $\lambda\leq\kappa$ and $X\subseteq
M^{P(\kappa^+)}$. We say that $X$ is morass-commutative beneath
$\lambda$ provided that for any $\zeta<\xi\leq\lambda$ and $g\in
\script{G}_{\zeta \xi}$, $x\in X\cap M^{P_{\zeta}}$ iff
$g(x)\in X\cap M^{P_{\xi}}$. We say that $X$ is morass-commutative if it is
morass-commutative beneath $\kappa$.
\end{definition}
\begin{definition}(Embedding-Commutative)
Let $\kappa$ be regular, $\gaptwo$ be a simplified
$(\kappa,2)$-morass, $\lambda\leq\kappa$ and
$X\subset M^P$.  We say that $X$ is
embedding-commutative beneath $\lambda$
provided that for any $\zeta<\xi\leq\lambda$ and
$f\in \script{F}_{\zeta \xi}$, $x\in X\cap M^{P_{\zeta}}$
iff $f_{\theta_{\zeta}}(x)\in X\cap M^{P_{\xi}}$.
We say that $X$ is embedding-commutative if it is
embedding-commutative beneath $\kappa$.
\end{definition}
Morass-commutativity and embedding-commutativity extend to relations and functions
on terms in the obvious way.
\begin{definition}(Grounded Order-Support)
Let $T_X\in M^P$ be forced to be a linear-ordering and
$X\subseteq M^P$ be a discerning set of terms for the
domain of $T_X$.  $X$ has grounded order-support
provided that for all $x, y\in X$ and
$G$, $P$-generic over $M$, there is $z\in M\cap X$ such that
$M[G]\models x<z<y$.
\end{definition}
If $X$ has grounded order-support, then for all $x, y\in X$,
$\forces \exists z\in M (x<z<y)$ and
for all $p\in P$ with $p\forces x<y$, there is
$z\in X\cap M$ and $q\leq p$ such that
$q\forces x<z<y$.  For instance, $\R$ has grounded support.
\begin{definition}(Upward Level-Dense)
Let $T_X\in M^P$ be forced to be a linear-ordering,
and $X\subseteq M^P$ be a set of discerning terms for
the domain of $T_X$.  $T_X$ is upward level-dense provide that
for every $x,y,z\in X$, in which
$z$ has strict order support $A\subseteq \kappa$, $p\in P$
with $p\forces x<z<y$ and $B\supseteq A$, there
is a discerning term $w$ with strict support $B$,
such that  $p\forces x<w<y$.
\end{definition}
If a set of discerning terms has grounded order-support, and terms bear
an order relation in a generic extension, then there
is an element of the ground model that is between the elements.  If a set of discerning
terms is upward level-dense, then between any pair
of elements there are elements of arbitrarily large
strict support (in the sense
of containment) between them.

We will require that the sequence of term functions
we are constructing
is closed under second components of morass-embeddings between
our ``fake morasses".  That is, for $\alpha<\beta\leq \omega_1$,
and $f\in \script{F}_{\alpha \beta}$, we require that
the partial term function constructed at
level $\alpha$, is closed under maps,
$f_{\theta_{\alpha}}:\varphi_{\theta_{\alpha}} \to \varphi_{\theta_{\beta}}$.

We revise the definition of morass-definability for extension
to the gap-2 construction.
\begin{definition}(Morass-Definable)\label{def4.10}
Let $\gaptwo$ be a simplified $(\omega_1,2)$-morass,
$P$ be the poset that adds generic reals
indexed by $\omega_3$,
$X\in M^P$ be a discerning
set of terms and $R\subset X\times X$.
We say $\langle X,R\rangle$ is a morass-definable $\eta_1$-ordering
(with respect to $\gaptwo$)
provided that
\begin{enumerate}
\item For every $A\subseteq\omega_2$,
$\langle X,R\rangle \cap (M^{P(A)})^3$ is
forced to be an $\eta_1$-ordering
\item $X$ and $R$ are morass-commutative and embedding-commutative
\item $X$ has grounded order-support and is upward level-dense
\item Every term of $X$ has countable support.
\end{enumerate}
If $G$ is $P$-generic over $M$ and
$\langle \bf{X},\bf{R}\rangle \in M[G]$ is an $\eta_1$-ordering
we say that $\langle \bf{X},\bf{R}\rangle \in M[G]$ is
morass-definable provided there is morass-definable
$\langle X,R\rangle$
with $val_G(\langle X,R\rangle)=\langle \bf{X},\bf{R}\rangle$.
\end{definition}
If $\langle X,R\rangle\cap{M^{P_{\theta_{\beta}}}}$
satisfies the conditions of Definition \ref{def4.10}
below $\beta$, then we say that $\langle X,Y\rangle$
is morass-definable below $\beta$.

If $M$ is a c.t.m. of $ZFC+CH$
containing a simplified $(\omega_1, 1)$-morass
and $M[G]$ is the Cohen extension adding $\aleph_2$
generic reals, then in $M[G]$, morass-definable
$\eta_1$-orderings are order-isomorphic.  Furthermore,
if $P$ is the poset for adding $\aleph_2$ generic reals,
and $(X,<_X)$ and $(Y,<_Y)$ are morass-definable and
are forced to be $\eta_1$-orderings, there is a level term
function, $\phi:X\to Y$, that is forced in all $P$-generic
extensions to be an order-isomorphism.
\section{The Gap-2 Construction}
In the proof that morass-definable $\eta_1$-orderings are
order-isomorphic in the Cohen extension adding $\aleph_2$-generic
reals, we built a term function with cardinality $\aleph_2$ by an
inductive construction of length $\omega_1$.  The constraint on
the length of the chain is governed by the possibility that
we are unable to to satisfy uncountably many simultaneous consistent order-constraints when
committing to discrete extensions.  The construction depends on
commutative extensions by morass maps to ``lift up" the
construction to exhaust the required commitments in the forcing
language adding $\aleph_2$ generic reals.  In the argument using the
gap-1 morass, all commitments can be met provided that any
countable subset of $\omega_2$ is in the range of a single morass
map from a countable vertex to the $\omega_1$ vertex.

This strategy does not extend to higher cardinality.
Instead, we work in a model of
$ZFC + CH$ that contains a simplified $(\omega_1,2)$-morass,
$\gaptwo$, using the morass-embeddings between ``fake" morasses,
that is, initial segments of the gap-two morass for vertices
$\theta_{\beta}$, $\beta \leq \omega_1$.
The maps $f_{\theta_{\beta}}:\varphi_{\theta_{\beta}} \to \omega_3$,
for embeddings $f\in \script{F}_{\beta \omega_1}$,
``lift" a term function on a domain of cardinality
$\aleph_1$ to a term function on a domain with cardinality
$\aleph_3$. Provided the term function constructed
at stage $\omega_1$ exhausts $X_{\omega_1}$ and $Y_{\omega_1}$, the
technical results of Section 3 allow us to complete
the construction of a term function on $X$ by closure under
second components of morass-embeddings of $\script{F}_{\beta \omega_1}$, $\beta<\omega_1$.

In Section 4 we adapted the key definitions from [\ref{Dumas}] to a higher
cardinality gap-1 morass and a gap-2 morass.  These definitions
concern technical considerations for building a function between
sets of terms in the forcing language for adding Cohen generic
reals that may be extended by compatible injections on
ordinals.  In passing from the gap-1 construction to the
gap-2 construction, we need to use the gap-1 morass maps
of $\script{G}$ and the second components of
morass-embeddings $\script{F}$, to extend
a partial construction of a term function below $\omega_1$ by
closure under second components of embeddings to a term function in the
forcing language adding generic reals indexed by $\omega_3$.
That is, we will make all discrete decisions
extending a function in an enhanced back-and-forth construction
of length $\omega_1$,
and ``lift" the entire construction by second components of embeddings of
$\script{F}_{\beta \omega_1}$ to exhaust terms with countable support in the
language adding generic reals indexed by $\omega_3$.

Assume that $\theta \leq \phi$ are ordinals, $f:\theta \to
\phi$ is an injection, and that $\tau$ is a term in the forcing
language adding generic reals indexed by $\theta$.  We define
$f(\tau)$ to be the term in the forcing language adding generic
reals indexed by $\phi$ that results from the formal substitution of
every indexing ordinal in $\theta$ appearing in $\tau$ by its image under $f$.
So $f(\tau)$ is a term in the forcing language adding
generic reals indexed by $f[\theta]\subseteq \phi$.  Morass maps of $\gapone$ and the
first and second components of morass
embeddings of $\gaptwo$, are order-preserving injections on
ordinals.  We treat morass maps and second components of
morass-embeddings as functions between terms of a forcing
language by this convention.

The morass maps of a simplified gap-1 morass have a restricted
character. There are only two morass maps from a morass vertex to its successor:
identity on an associated ordinal, and a single ``splitting" map
that translates the cofinal end segment of that ordinal.
We need to consider a new collection of ordinal
maps along the morass that will enrich the terms we
can define by closure under ordinal injections, the second components of the
morass embeddings, $f_{\theta_{\alpha}}:\varphi_{\theta_{\alpha}}\to \varphi_{\theta_{\beta}}$, for
$f\in \script{F}_{\alpha \beta}$ ($\alpha<\beta\leq\omega_1)$.
We are able to restrict
our attention to the morass-embeddings of the amalgamations
$\script{F}_{\alpha \alpha+1}$.
Each $f_{\theta_{\alpha}}$ is an
order-preserving injection, but may split its domain into many
disconnected pieces that have less predictable intersections with
morass maps.
\begin{definition}(Embedding-closure)
Let $\beta<\gamma \leq \omega_1$
and $X$ be a set of terms in the forcing language $M^P_{\gamma}$.
The closure of $X$ under $\script{F}_{\beta \gamma}$ is the smallest
set of terms of $M^P_{\gamma}$ such that
$x\in X\cap M^{P_{\beta}}$ and $f\in \script{F}_{\beta \gamma}$ implies
$f_{\theta_{\beta}}(x)\in X$.  The set $X$ is closed under $\script{F}_{\beta \gamma}$ if
the closure of $X$ is $X$.  If $\beta<\omega_1$,
and $X$ is closed under $\script{F}_{\beta \beta+1}$, then we say that $X$ is
embedding-closed at level $\beta$.
\end{definition}
We will require that our sequence of functions
be embedding-closed at level $\beta$ for all $\beta<\omega_1$.  We prove a lemma that
allows for a transfinite construction of a sequence of embedding-closed
term functions
along a simplified gap-2 morass.
If $X$ is a morass-definable $\eta_1$-ordering,
Then $X$ has grounded order-support.
If $D\subset X$ is countable, then c.c.c
of $P$ implies that there is a countable
extension of $D$ by elements of the ground model
that has grounded order-support.
\begin{lemma}\label{lem5.3}
Let $\beta$ be a countable ordinal and assume
\begin{enumerate}
\item $\gaptwo$ is a simplified $(\omega_1,2)$-morass
\item $X$ and $Y$ are morass-definable $\eta_1$-orderings
\item $D\subseteq X\cap M^{P_{\theta_{\beta}}}$
has grounded order-support and is morass-commutative
\item $F:D\to Y$ is a morass-commutative level term injection that is
forced to be order-preserving in all generic extensions.
\end{enumerate}
Then $\phi=\bigcup \{ f_{\theta_{\beta}}[F]\mid f\in \script{F}_{\beta \beta+1} \}$
is a level term injection that is forced to be an order-preserving injection.
\end{lemma}
Proof:
Let $\gaptwo$, $X$, $Y$, $D$ and $F$ satisfy the
hypotheses of the lemma, and $\phi=\bigcup \{ f_{\theta_{\beta}}[F]\mid f\in
\script{F}_{\beta \beta+1} \} \subseteq X\times Y$.
The family of embeddings, $\script{F}_{\beta \beta+1}$
is an amalgamation, so it is composed of a single right-branching $f\in \script{F}_{\beta \beta+1}$,
and all possible left-branching embeddings.  For each
left branching embedding of $g\in \script{F}_{\beta \beta+1}$,
$g_{\theta_{\beta}}\in \script{G}_{\theta_{\beta} \theta_{\beta+1}}$.
The second components of embeddings of $\script{F}$
are injections.  Therefore $\phi$ is level, and the elements
of the range of $\phi$ are discerning terms in $Y$.

To see that $\phi$ is forced to be a well-defined function,
suppose that\\
$(x_1,y_1), (x_2, y_2) \in F$; $f,g\in\script{F}_{\beta \beta+1}$;
$p\in M^{P_{\theta_{\beta+1}}}$
and $p\forces f_{\theta_{\beta}}(x_1)=g_{\theta_{\beta}}(x_2)$.  Since $x_1$ and
$x_2$ are discerning terms, so are $f_{\theta_{\beta}}(x_1)$
and $g_{\theta_{\beta}}(x_2)$.  Therefore $f_{\theta_{\beta}}(x_1)$
and $g_{\theta_{\beta}}(x_2)$ must have identical strict support, $S$,
and $p$ may be replaced by a condition of $Fn(S\times \omega,2)$.
We claim that $p\forces f_{\theta_{\beta}}(y_1)=g_{\theta_{\beta}}(y_2)$.
We observe that for $\alpha\in S$, there are $\alpha_1, \alpha_2\in \varphi_{\theta_{\beta}}$
such that $f_{\theta_{\beta}}(\alpha_1)=g_{\theta_{\beta}}(\alpha_2)=\alpha$.
Then there is $\bar{\alpha}\in \theta_{\eta}$ and
$\bar{h}\in \script{G}_{\eta \theta_{\beta}}$ such that
$\bar{h}(\bar{\alpha})=\alpha_1$.
However,
\[ f_{\theta_{\beta}}\comp \bar{h}(\bar{\alpha})=f_{\eta \theta_{\beta}}(\bar{h})\comp f_{\eta}(\bar{\alpha}). \]
So $f_{\eta}(\bar{\alpha})=\alpha_2$ and $f_{\eta \theta_{\beta}}(\bar{h})(\alpha_2)=g_{\theta_{\beta}}(\alpha_2)$.
By Lemma \ref{lem3.2}, $f_{\eta \theta_{\beta}}(\bar{h})$ and $g_{\theta_{\beta}}$ are equal up to $\alpha_2$.
Specifically, if $\gamma\in \varphi_{\eta}$, $h\in \script{G}_{\eta \theta_{\beta}}$ and $h(\gamma)\in supp(x_1)$,
then $f_{\theta_{\beta}}\comp h(\gamma)=g_{\theta_{\beta}}\comp f_{\eta}(\gamma)$.

Define a relation $G\subseteq \varphi_{\eta}\times \varphi_{\theta_{\beta}}$ by $(\gamma,\lambda)\in G$
iff there is $h\in \script{G}_{\eta \theta_{\beta}}$ with $\lambda\in supp(x_1)$ and $h(\gamma)=\lambda$.
We claim that $G$ is a well-defined function and order-preserving injection.  If $h_1$ and $h_2$ are
witnessing functions from the definition of $G$ for $\gamma$, then
$f_{\theta_{\beta}}\comp h_1(\gamma)=g_{\theta_{\beta}}\comp f_{\eta}(\gamma)=f_{\theta_{\beta}}\comp h_2(\gamma)$.
Since $f_{\theta_{\beta}}$ is an injection, $h_1(\gamma)=h_2(\gamma)$, and $G$ is a well-defined function.
Assume that $\gamma_1, \gamma_2\in \varphi_{\eta}$ and $G(\gamma_1)=G(\gamma_2)=\lambda$.
Then
\[ f_{\theta_{\beta}}(\lambda)=g_{\theta_{\beta}}\comp f_{\eta}(\gamma_1)=
g_{\theta_{\beta}}\comp f_{\eta}(\gamma_2). \]
Therefore $\gamma_1=\gamma_2$ and $G$ is an injection.
So
\[ f_{\theta_{\beta}}((x_1,y_1))=(f_{\theta_{\beta}}\comp G)\comp G^{-1}((x_1,y_1))=
g_{\theta_{\beta}}\comp f_{\eta}((G^{-1}(x_1),G^{-1}(y_1))). \]
Therefore $f_{\eta}\comp G^{-1}(x_1)=x_2$.
Since $F$ is morass-commutative, and $D$ has grounded order-support,
$\phi \restrict_{X_{\eta}}$ is morass-commutative.  Hence
\[ f_{\theta_{\beta}}(y_1)=g_{\theta_{\beta}}(y_2). \]
Therefore $\phi$ is forced to be a well-defined function.  The proof that $\phi$ is
forced to be an injection is essentially identical.

We show that $\phi$ is forced to be order-preserving.
Consider $\phi$ as a function on $\{ f_{\theta_{\beta}}[D]\mid f\in \script{F}_{\beta \beta+1} \}$.
Let $f, g\in \script{F}_{\beta \beta+1}$.
If $x\in M\cap D$, then $\phi(x)=F(x)$.
Assume $x_1, x_2\in D$ and $p\forces x_1<x_2$.
Let $y_1=F(x_1)$ and $y_2=F(x_2)$ and
Since $D$ has a grounded order-support,
there is $z$ in $M\cap D$ and $q\leq p$ such that
$q\forces x_1<z<x_2$.
Then
\[ q\forces f_{\theta_{\beta}}(x_1)<z<g_{\theta_{\beta}}(x_2). \]
So
\[ q\forces \phi(f_{\theta_{\beta}}(x_1))<\phi(z)=F(z)<\phi(g_{\theta_{\beta}}(x_2)). \]
Therefore it is forced that $\phi$ is order-preserving.
\halmos
\begin{theorem}\label{thm5.3}
Let $M$ be a model of $ZFC+CH$ containing a
simplified $(\omega_1,2)$-morass, $\gaptwo$, and
$P=Fn(\omega_3\times \omega, 2)$.
Let $X$ and $Y$ be sets of discerning terms
in $M^{P}$ for
morass-definable $\eta_1$-orderings
(with respect to $\gaptwo$).
Then there is a level function from
$X$ to $Y$ that is forced to be an
order-isomorphism.
\end{theorem}
Proof.
For $\alpha<\omega_1$,
let $X_{\alpha}=X\cap M^{P_{\alpha}}$
and $Y_{\alpha}=Y\cap M^{P_{\alpha}}$.
We consider $X_{\alpha}$ and $Y_{\alpha}$ as the restrictions of
$X$ and $Y$, resp., to the forcing
language adding generic reals indexed by $\varphi_{\theta_{\alpha}}$.
In any $P$-generic extension of $M$, $M[G]$,
the interpretation of $X_{\alpha}$ in $M[G]$ is
the interpretation of $X_{\alpha}$ in $M[G_{\alpha}]$
where $G_{\alpha}$ is the factor of $G$ that is
$P_{\alpha}$-generic over $M$.

We will use the previous Lemma to construct a
morass-commutative level term bijection from $X_{\omega_1}$
to $Y_{\omega_1}$ that is forced
to be an order-isomorphism.  The closure
under embeddings, $f_{\theta_{\beta}}$ where $f\in \script{F}_{\beta \omega_1}$, of this
function will be the term function we seek.

Let $\langle x_{\beta} \mid \beta<\omega_1 \rangle$
be an enumeration of $X_{\omega_1}$ that satisfies the condition
$x_{\alpha}\in X_{\alpha}$ for all $\alpha<\beta$.
Let $\langle y_{\beta} \mid \beta<\omega_1 \rangle$ satisfy the same condition
with respect to $Y_{\omega_1}$.

We will inductively construct a transfinite sequence of
morass-commutative term functions $\langle F_{\beta}:D_{\beta}\to E_{\beta}\mid \beta<\omega_1 \rangle$
that satisfies the following for all $\alpha \leq \beta<\omega_1$,
\begin{enumerate}
\item $D_{\beta}\subseteq X_{\theta_{\beta}}$ and $E_{\beta}\subseteq Y_{\theta_{\beta}}$
are morass-commutative, countable sets of terms with grounded order support
\item $D_{\alpha}\subseteq D_{\beta}$ and $E_{\alpha}\subseteq E_{\beta}$
\item $x_{\beta}\in D_{\beta}$ and $y_{\beta}\in E_{\beta}$
\item $F_{\beta}$ is a morass-commutative level term function that is forced to be an order-preserving bijection
\item $f_{\theta_{\alpha}}[F_{\alpha}]\subseteq F_{\beta}$ for all
$f\in \script{F}_{\alpha \beta}$
\end{enumerate}
We call a sequence of term functions satisfying these conditions (beneath $\beta$)
and extendable sequence.  We argue be induction on $\gamma<\omega_1$

Case 1: $\gamma=0$.  Then $x_0$ and $y_0$ are elements of the ground model, $M$.
Let $F_0=\{ (x_0,y_0)\}$.

Case 2: $\gamma=\beta+1$.  Let $\langle F_{\alpha}:D_{\alpha}\to E_{\alpha}
\mid \alpha \leq \beta \rangle$
be an extendable sequence.
By Lemma \ref{lem5.3}, $(\bigcup_{f\in \script{F}_{\beta \gamma}} f_{\theta_{\beta}}[F_{\beta}]):
\bigcup_{f\in \script{F}_{\beta \gamma}} f_{\theta_{\beta}}[D_{\beta}]
\to \bigcup_{f\in \script{F}_{\beta \gamma}} f_{\theta_{\beta}}[E_{\beta}]$ is forced
to be an order-preserving bijection.
Let $S=supp(x_{\gamma})$ and $\alpha \leq \gamma$ be the least ordinal for which
there exists $\sigma \in G_{\alpha \varphi_{\gamma}}$ and $\bar{S}\subseteq \varphi_{\alpha}$
such that $\sigma[\bar{S}]=S$.  Then there is $\bar{x}\in M^{P(\varphi_{\alpha})}$
such that $\sigma(\bar{x})=x_{\gamma}$.
Let $D^*=\{ \sigma(\bar{x}) \mid \sigma\in \script{G}_{\alpha \theta_{\gamma}}) \}\subseteq X$.
Then $D^*$ is the smallest morass-commutative subset of $X_{\beta+1}$ that contains $x_{\gamma}$.
Let $D'\subset X$
be a countable subset of $X_0$
so that $D=D'\union D^*$ has grounded order-support.
By repeated applications of Lemma 4.5 of [\ref{Dumas}], there is
$F':(D'\cup D_{\beta}) \to Y$ that is a morass-commutative  level term function that
is forced to be order-preserving and $D'\cup D_{\beta}$ has grounded order-support.  By Lemma 4.5, [\ref{Dumas}], there is
$\bar{y}\in Y$ such that $F'\cup \{ (\bar{x},\bar{y}) \}$ is a level term function
that is forced to be order-preserving.
Let $F=F'\cup \{ f_{\theta_{\beta}}\comp g((\bar{x},\bar{y})) \mid
f\in \script{F}_{\beta \gamma}, g\in \script{G}_{\alpha \theta_{\beta}} \}$.
Then $F$ is a level term function that is forced to be a bijection.  To see that
$F$ is order-preserving, assume that
$f,h\in \script{F}_{\beta \gamma}$,
$g_1,g_2\in \script{G}_{\alpha \theta_{\beta}}$,
$x_1=f_{\theta_{\beta}}\comp g_1(\bar{x})$,
$x_2=h_{\theta_{\beta}}\comp g_2(\bar{x})$,
$G$ is $P_{\gamma}$-generic over $M$ and $M[G]\models x_1<x_2$.
Since the domain of $F$ has grounded order-support, there is $x_0\in D'$ such that
$M[G]\models x_1<x_0<x_2$.  Then
\[ M[G]\models F(x_1)<F(x_0)<F(x_2). \]
Therefore $F$ is forced to be order-preserving.

Let $E$ be the range of $F$.
Let $E^*$ be the smallest morass-commutative subset of $Y_{\beta+1}$ that contains $y_{\gamma}$.
Then $E$ has grounded order-support and there is a countable extension
of $E \cup E^*$ by elements of $Y_0$, $E_{\gamma}$,
so that $E_{\gamma}$ has grounded order-support.
Again, by applications of Lemma 4.5 [\ref{Dumas}], there is a level injection, $F_{\gamma}$, such that
$F^{-1}_{\gamma}:E_{\gamma}\to X$ is a level injection extending
$F^{-1}$ that is forced to be an order-preserving injection.
Let $D_{\gamma}$ be the range of $F^{-1}_{\gamma}$.  Then
$\langle F_{\alpha}\mid \alpha \leq \gamma \rangle$ is an extendable sequence.

Case 3: $\gamma$ a limit ordinal.  It is routine to verify that
$F=\bigcup_{\alpha<\gamma, f\in \script{F}_{\alpha \gamma}} f_{\theta_{\alpha}}[F_{\alpha}]$
is a level injection that is forced to be an order-preserving injection.
In a manner identical to the successor case, $F$ may be extended to a
level injection, with domain containing $x_{\gamma}$ and range
containing $y_{\gamma}$, both bearing grounded order-support, in which $F$ is forced
to be an order-preserving injection.

Let $F=\bigcup_{\alpha<\omega_1, f\in \script{F}_{\alpha \omega_1}} f_{\theta_{\alpha}}[F_{\alpha}]$.
Then $F:X\to Y$ is a level bijection that is forced to be an order-isomorphism.
By Lemma \ref{lem3.9}, the domain of $F$ is $X$ and the range of $F$ is $Y$.
\halmos
\section{Ultrapowers of $\R$ over $\omega$}
We turn our attention to $\romu$, an ultrapower of $\R$ over a
non-principal ultrafilter on $\omega$, $U$. By results of G. Dales
[\ref{Dales1}], J. Esterle [\ref{Esterle2}] and B. Johnson
[\ref{Johnson}] the existence of an $\R$-linear order-preserving
monomorphism from the finite elements of $\romu$, for some
non-principal ultrafilter on $\omega$, $U$, into the Esterle
algebra is sufficient to prove the existence of a discontinuous
homomorphism of $C(X)$, the algebra of continuous real valued
functions on $X$, where $X$ is an infinite compact Hausdorff
space.  In a later paper, relying on the techniques of
this paper, we prove that it is consistent that such a
monomorphism exists in a model of set theory in which
$2^{\aleph_0}=\aleph_3$.  In anticipation of such a construction,
we finish this paper with a proof that in the Cohen extension
adding $\aleph_3$-generic reals of a model of ZFC+CH containing a
simplified $(\omega_1,2)$-morass, there is an ultrapower of $\R$ over
a standard ultrafilter on $\omega$ that is gap-2 morass-definable, and
hence is order-isomorphic with other gap-2 morass-definable
$\eta_1$-orderings. This result extends Theorem 6.15
[\ref{Dumas}], that in the Cohen extension adding
$\aleph_2$-generic reals of a model of $ZFC + CH$ containing a
simplified gap-1 morass, ultrapowers of $\R$ over standard
ultrafilters on $\omega$ are order-isomorphic.

By Theorem \ref{thm5.3}, in order to prove that an ultrapower of
$\R$, $\romu$, is order-isomorphic with a gap-2 morass-definable
$\eta_1$-ordering, it sufficient to prove that $\romu$ is gap-2
morass-definable. We show that this is so, provided that $U$ is a
standard non-principal ultrafilter. We adapt the definition of
standard ultrafilter (Definition 6.14 [\ref{Dumas}]) to a
simplified gap-2 morass.
\begin{definition}(Standard Term for a subset of $\omega$)
A standard term for a subset of $\omega$ is a term $x\in M^P$,
such that for each $(\tau, p) \in x$, $\tau$ is a canonical term
in $M^P$ for a natural number.
\end{definition}
\begin{definition}(Standard Term for an Ultrafilter)
Let $\lambda \leq \omega_3$ and $U\in M^{P({\lambda})}$ be a
morass-commutative and embedding-commutative
set of standard terms for subsets of $\omega$ (below
$\lambda$) such that for all $S\subseteq \lambda$, $U\cap
M^{P(S)}$ is forced to be an ultrafilter in all
$P(S)$-generic extensions for $M$.  Then $U$ is a standard
term for an ultrafilter below $\lambda$.
\end{definition}
If $U$ is a standard term for an ultrafilter below $\omega_3$, we
say it is a standard ultrafilter.
\begin{definition}(Complete Standard Term for an Ultrafilter)
Let $U$ be a standard term for an ultrafilter below $\lambda$. $U$
is complete provided that for every standard term for a subset of
$\omega$, $u\in M^{P(\lambda)}$, $u\in U$ iff $\forces u\in U$.
\end{definition}
That is, below $\lambda$, every standard term for a subset of
$\omega$ that is forced to be in $U$, is a member of $U$. Every
standard term for an ultrafilter has a complete extension. Let $U$
be a complete standard term for an ultrafilter and $u$ be a
standard term for a subset of $\omega$. There is a standard term
for a subset of $\omega$ that decides the membership of $u$ in $U$
in all $P(\lambda)$-generic extensions.  Since $U$ is forced to be
an ultrafilter, $\{ p\in P(\lambda) \mid p\forces u\in U \vee
p\forces u\nin U\}$ is dense in $P(\lambda)$.  We consider $u$ as
a term for a binary sequence.  In this sense, $p\forces n\in u$
iff $p\forces u_n=1$. It is clear that there is a standard term
for a binary sequence, $v$, such that $p\forces u_n=1$ iff
$p\forces v_n=0$.  Let $d(u)$ be the standard term for a subset of
$\omega$ so that it is forced that $d(u)=u$ if $p\forces u\in U$
and $d(u)=v$ if $p\forces u\nin U$.  Then it is forced in all
generic extensions that $d(u)\in U$.

We wish to construct a standard term for an ultrafilter that
commutes with the second components of embeddings of $\script{F}_{\beta \beta+1}$,
for all $\beta<\omega_1$.
That is, if $\gaptwo$ is a simplified gap-2 morass, we
construct a sequence of standard terms for ultrafilters, $\langle U_{\beta}\mid\beta<\omega_1\rangle$,
where for each successor, $\beta$, $U_{\beta}$ is a set of standard
terms for subsets of $\omega$ in the language adding generic
reals indexed by $\varphi_{\theta_{\beta}}$, that is forced to be a
non-principal ultrafilter
in all $P_{\beta}$-generic extensions of $M$.  Furthermore,
we require that for all countable $\beta$, and all $f\in \script{F}_{\beta \beta+1}$,
$f_{\theta_{\beta}}[U_{\beta}]\subseteq U_{\beta+1}$.  Since
$\script{F}_{\beta \beta+1}$ is an amalgamation, all $f\in \script{F}_{\beta \beta+1}$,
with a single exception, are left-branching, and hence
$f_{\theta_{\beta}}\in \script{G}_{\theta_{\beta} \theta_{\beta+1}}$.
By results of [\ref{Dumas}], the morass closure of $U_{\theta_{\beta}}$
under left-branching embeddings have f.i.p.  The single right-branching
embedding of $\script{F}_{\beta \beta+1}$ must be handled independently.
\begin{theorem}
There is a standard ultrafilter, $U$, that
commutes with the embeddings of $\script{F}_{\beta \beta+1}$,
for all $\beta<\omega_1$.
The ultrapower, $\romu$, is a gap-2 morass-definable
$\eta_1$-ordering.
\end{theorem}
Proof.  By Lemma \ref{lem3.12} the morass-embeddings of an
amalgamation $\script{F}_{\beta \beta+1}$, for $\beta<\omega_1$, are compatible.
We construct the sequence $\langle U_{\beta} \rangle$ by induction on $\beta$.
The limit case is routine, so we assume that $\langle U_{\alpha}\mid \alpha \leq \beta \rangle$ has been defined.
We wish to show that the closure of $U_{\beta}$ under the embeddings of $\script{F}_{\beta \beta+1}$
has f.i.p.
\begin{lemma}\label{lem6.5}
Let $U_{\beta}$ be a standard ultrafilter, and $\script{F}_{\beta \beta+1}$ be an
amalgamation.  Let $f\in \script{F}_{\beta \beta+1}$ be the right-branching embedding,
and $g_1,\ldots,g_n \in \script{G}_{\theta_{\beta} \theta_{\beta+1}}$.
If $u_1,\ldots,u_n,v\in U_{\beta}$,
then $\forces g_1(u_1)\cap \ldots g_n(u_n)\cap f_{\theta_{\beta}}(v)\neq \emptyset$.
\end{lemma}
Proof of Lemma.
Let $u$ be a standard term for a subset of $\omega$ such that $\forces u=u_1 \cap \ldots \cap u_n \cap v$.
Then $\forces u\in U_{\theta_{\beta}}$.
If $U_{\theta_{\beta}}$ is complete, then $u\in U_{\theta_{\beta}}$.
If $k\in \N$, $p\in P_{\beta}$ and $h\in \script{F}_{\beta \beta+1}$,
then $p\forces k\in u$ iff $h_{\theta_{\beta}}(p)\forces k\in h_{\theta_{\beta}}(u)$.
Hence if $p\forces k\in u$,
then for all $i\leq n$,
\[  g_i(p)\forces k\in g_i(u). \]
By Lemma \ref{lem3.12}, $f_{\theta_{\beta}}(p), g_1(p),\ldots,g_n(p)$
are compatible, as members of $P$.
It follows that
\[ f_{\theta_{\beta}}(p)\cdot g_1(p)\cdots g_n(p)
\forces k\in g_1(u)\cap \dots \cap g_n(u)\cap f_{\theta_{\beta}}(u). \]
If there were a generic extension in which
$g_1(u_1)\cap \ldots \cap g_n(u_n)\cap f_{\theta_{\beta}}(v)=\emptyset$, then
in that generic extension,
\[ g_1(u)\cap \dots \cap g_n(u)\cap f_{\theta_{\beta}}(u)=\emptyset. \]
Suppose a condition, $p\in P_{\beta+1}$, forced this.  Let $S\subseteq f_{\theta_{\beta}}[\varphi_{\theta_{\beta}}]$
be the support of $f_{\theta_{\beta}}(u)$.  If $\pi$ were
the projection onto $S$, then by a straightforward generalization
of Lemma 6.2 [\ref{Dumas}],
\[ \pi(p)\forces \pi(g_1(u)\cap \ldots \cap g_n(u))\cap \pi(f_{\theta_{\beta}}(u))=\emptyset. \]
Then,
\[ \pi(p)\forces g_1(u)\cap \ldots \cap g_n(u)\cap f_{\theta_{\beta}}(u)=\emptyset. \]
However there is a condition $q\leq f^{-1}_{\theta_{\beta}}(\pi(p))$ and $k\in \N$
such that $q\forces k\in u$.
Therefore there is a condition $r\leq g_1(q)\cdots g_n(q)\cdot f_{\theta_{\beta}}(q)\leq \pi(p)$
where
\[ r\forces k\in g_1(u)\cap \ldots \cap g_n(u)\cap f_{\theta_{\beta}}(u). \]
So the intersection of $g_1(u_1), \ldots ,g_n(u_n), f_{\theta_{\beta}}(v)$
is forced in all generic extensions to be nonempty.  This completes the proof of the Lemma.
\halmos
Continuing the proof of the Theorem,
it follows from Lemma \ref{lem6.5} that the union under the embeddings of $\script{F}_{\beta \beta+1}$
of a standard ultrafilter in $M^{P_{\beta}}$ has f.i.p., and may be extended to
a standard term for an ultrafilter in the forcing language adding
generic reals indexed by $\varphi_{\theta_{\beta+1}}$.  It is
routine to see that there is a standard ultrafilter, $U\subseteq M^{P_{\omega_2}}$, that commutes
with the embeddings of
$f\in \script{F}_{\alpha \beta}$ for $\alpha<\beta\leq \omega_1$.
We claim that $\romu$ is a gap-2 morass-definable $\eta_1$-ordering.
It was shown in $[\ref{Dumas}]$ that $\romu$ is
upward level-dense, has countable support
and is morass-commutative.  By an application of the proof
of Lemma 6.10 [\ref{Dumas}], $\romu$ is embedding-commutative.
We show that $\romu$ has grounded order support.
We repeat the proof of Lemma 6.13 [\ref{Dumas}].
Let $x$ and $y$ be terms for sequences of reals, $p\in P$
and $p\forces [x]_U<[y]_U$.
There are terms for sequences or reals, $\bar{x}$ and $\bar{y}$,
such that $p\forces [x]_U=[\bar{x}]_U$, $p\forces [y]_U=[\bar{y}]_U$
and $p\forces \all i\in \N (\bar{x}_i<\bar{y}_i)$.
Let $\theta$ be an ordinal and $\sigma:\varphi_{\theta_{\beta+1}} \to \theta$
be an injection,
where $\varphi_{\theta_{\beta+1}}\cap \theta=\emptyset$.
Then, since $\R$ is level dense,
\[ p\cdot \sigma(p)\forces \forall i\in \N (\bar{x}_i<\sigma(\bar{y}_i)). \]
Let $a_i=sup\{ q\in \Q \mid (\exists r<p)\wedge (r\forces q<\bar{x}_i) \}$
and $b_i=inf\{ q\in \Q \mid (\exists r<p) (\sigma(r)\forces \sigma(y_i)<q) \}$.
Then
\[ \forall i\in \N (a_i\leq b_i). \]
For $i\in \N$, let $c\in \ro \cap M$ be such that
$a_i\leq c_i \leq b_i$. Let $U'$ be a standard ultrafilter extending $U\cup \sigma(U)$.
Then
\[ p\cdot \sigma(p) \forces [x]_{U'}=[\bar{x}]_{U'}\leq [c]_{U'} \leq [\sigma(\bar{y})]_{U'}=[\sigma(y)]_{U'}. \]
In the language of the gap-1 argument, $\romu$ has a grounded order base.
Therefore, by Lemma 4.4 [\ref{Dumas}], $\romu$ has a grounded order support.
Hence, $\romu$ is gap-2 morass-definable.
\halmos
\section{Next Results}
In this paper and [\ref{Dumas}] we have extended the classical
result that $\eta_1$-orderings of cardinality $\aleph_1$ are
order-isomorphic.  In the next paper we extend the result that
there is an $\R$-linear isomorphism between $\eta_1$-ordered
real-closed fields of cardinalty $\aleph_1$. We show that in the
Cohen extension adding $\aleph_2$-generic reals to a model of
ZFC+CH containing a simplified $(\omega_1,1)$-morass, there is a
morass-definable $\R$-linear isomorphism between $\eta_1$-ordered
morass-definable real-closed fields. We then show that in the
Cohen extension adding $\aleph_3$-generic reals to a model of
ZFC+CH containing a simplified $(\omega_1,2)$-morass, there is a
gap-2 morass-definable $\R$-linear isomorphism between gap-2
morass-definable $\eta_1$-orderings. With these results we
are able to extend the theorem of Woodin [\ref{Woodin}] that it is
consistent that $2^{\aleph_0}=\aleph_2$ and there exists a
discontinuous homomorphism of $C(X)$. We show that it is
consistent that $2^{\aleph_0}=\aleph_3$ and that there exists a
discontinuous homomorphism of $C(X)$, for any infinite compact
Hausdorff space, $X$.

\end{document}